\documentclass[11pt]{article}

\usepackage{xspace}
\usepackage{url}
\usepackage{mathtools}
\usepackage{amssymb}
\usepackage{amsthm}
\usepackage{empheq}
\usepackage{latexsym}
\usepackage{enumitem}
\usepackage{graphicx}
\usepackage{eurosym}
\usepackage{dsfont}
\usepackage{appendix}
\usepackage{color} 
\usepackage[unicode]{hyperref}
\usepackage{frcursive}
\usepackage[utf8]{inputenc}
\usepackage[T1]{fontenc}
\usepackage{geometry}
\usepackage{multirow}
\usepackage{todonotes}
\usepackage{lmodern}
\usepackage{anyfontsize} 
\usepackage{stmaryrd}
\usepackage{cleveref}
\usepackage[english]{babel}
\usepackage[english=british]{csquotes}
\usepackage{mathtools}
\usepackage{accents}
\usepackage{color}
\usepackage{hyperref}
\usepackage{bm}
\usepackage{subcaption}
\usepackage{algorithm}
\usepackage{algpseudocode}

\definecolor{red}{rgb}{0.7,0.15,0.15}
\definecolor{green}{rgb}{0,0.5,0}
\definecolor{blue}{rgb}{0,0,0.7}
\hypersetup{colorlinks, linkcolor={red}, citecolor={green}, urlcolor={blue}}

\allowdisplaybreaks

\usepackage{natbib}
\setcitestyle{numbers,open={[},close={]}}

\definecolor{red}{rgb}{0.7,0.15,0.15}
\definecolor{green}{rgb}{0,0.5,0}
\definecolor{blue}{rgb}{0,0,0.7}
\hypersetup{colorlinks, linkcolor={red}, citecolor={green}, urlcolor={blue}}
			
\makeatletter \@addtoreset{equation}{section}

\newtheorem{theorem}{Theorem}[section]
\newtheorem{assumption}[theorem]{Assumption}

\newtheorem{lemma}[theorem]{Lemma}
\newtheorem{proposition}[theorem]{Proposition}

\newtheorem{remark}[theorem]{Remark}

\DeclareUnicodeCharacter{014D}{\=o}
\newcommand{\ba}{\begin{array}}
\newcommand{\ea}{\end{array}}
\newcommand{\be}{\begin{equation}}
\newcommand{\ee}{\end{equation}}
\newcommand{\bea}{\begin{eqnarray}}
\newcommand{\eea}{\end{eqnarray}}
\newcommand{\beaa}{\begin{eqnarray*}}
\newcommand{\eeaa}{\end{eqnarray*}}

\def\dbE{\mathbb{E}}
\def\dbF{\mathbb{F}}

\def\dbN{\mathbb{N}}
\def\dbP{\mathbb{P}}
\def\dbR{\mathbb{R}}

\def\a{\alpha}
\def\b{\beta}
\def\g{\gamma}
\def\d{\delta}
\def\e{\varepsilon}

\def\l{\lambda}
\def\m{\mu}

\def\si{\sigma}
\def\t{\tau}
\def\f{\varphi}
\def\th{\theta}

\def\D{\Delta}

\def\cF{{\cal F}}

\def\cJ{{\cal J}}

\def\cN{{\cal N}}
\def\cO{{\cal O}}

\def\cT{{\cal T}}

\def\q{\quad}

\def\qed{ \hfill \vrule width.25cm height.25cm depth0cm\smallskip}

\def\bx{{\boldsymbol{x}}}

\def\bX{{\boldsymbol{X}}}

\def\1{\mathbf{1}}

\def\bI{{\boldsymbol{I}}}
\def\bi{{\boldsymbol{i}}}
\def\bI{{\boldsymbol{I}}}
\def\bi{{\boldsymbol{i}}}

\begin{document}

\title{Deep Learning for the Multiple Optimal Stopping Problem}
\author{Mathieu Laurière\footnote{Shanghai Center for Data Science; NYU-ECNU Institute of Mathematical Sciences at NYU Shanghai; NYU Shanghai, Shanghai,
People’s Republic of China, mathieu.lauriere@nyu.edu} \quad Mehdi Talbi\footnote{Laboratoire de Probabilités, Statistiques et Modélisation, Université Paris-Cité, France, mtalbi@lpsm.paris}}

\date{}

\maketitle

\begin{abstract}
This paper presents a novel deep learning framework for solving multiple optimal stopping problems in high dimensions. While deep learning has recently shown promise for single stopping problems, the multiple exercise case involves complex recursive dependencies that remain challenging. We address this by combining the Dynamic Programming Principle with neural network approximation of the value function. Unlike policy-search methods, our algorithm explicitly learns the value surface. We first consider the discrete-time problem and analyze neural network training error. We then turn to continuous problems and analyze the additional error due to the discretization of the underlying stochastic processes. Numerical experiments on high-dimensional American basket options and nonlinear utility maximization demonstrate that our method provides an efficient and scalable method for the multiple optimal stopping problem.
\end{abstract}

\section{Introduction}

\noindent\textbf{Classical Optimal Stopping.} Optimal stopping problems consist of choosing a stopping time for a stochastic process to maximize a specific payoff, which is defined as a function of the process at the chosen time. The decision must be non-anticipative; that is, the stopping time must be adapted to the filtration generated by the underlying process. This framework models a diverse range of decision-making processes, including the pricing of American derivatives in finance, the evaluation of investment opportunities in economics (real options), and sequential hypothesis testing in statistics. 
It is well known that the (first) optimal stopping time is given by the first time the reward process reaches its Snell envelope, defined as the smallest supermartingale dominating the reward process (see e.g.\ \citeauthor*{shiryaev2008optimal} \cite{shiryaev2008optimal} or \citeauthor*{karatzas1998methods} \cite{karatzas1998methods} for a general overview). This characterization stems from the Dynamic Programming Principle (DPP), which underpins most numerical approaches. Notable classical examples include the celebrated Longstaff-Schwarz method \cite{longstaff2001valuing}, based on recursive least squares regressions, and the stochastic mesh approach of \citeauthor*{broadie2004stochastic} \cite{broadie2004stochastic}. Extensions of these methods include the use of statistical learning theory~\cite{egloff2005monte} and path signatures for non-Markovian settings~\cite{bayer2023optimal,bayer2025primal}.

\noindent\textbf{Deep Learning for Optimal Stopping.}
Recently, deep learning has emerged as a powerful tool to overcome the curse of dimensionality in optimal stopping. \citeauthor*{becker2019deep} \cite{becker2019deep} proposed a seminal approach combining dynamic programming with a neural network approximation of the stopping rule. This was later generalized to higher dimensions and value function computation in~\cite{becker2021solving}, and analyzed theoretically in~\cite{gonon2024deep} and~\cite{bayer2021randomized} regarding convergence rates and randomized stopping.
Alternative deep learning formulations have also gained traction. These include approaches based on Backward Stochastic Differential Equations (BSDEs)~\cite{gao2023convergence,yang2024deep,wei2024deep}, methods that learn the stopping boundary directly~\cite{reppen2025neural,soner2025stopping}, and primal-dual bounds via deep learning~\cite{guo2025simultaneous}. Beyond standard neural networks, literature explores Gaussian processes~\cite{dwarakanath2022optimal,ludkovski2025optimal}, Deep Q-Learning variants~\cite{ery2021solving,fathan2021deep}, and recurrent neural networks for non-Markovian problems~\cite{damera2023deep}. These techniques have found broad application in finance (see e.g.~\cite{hu2020deep,felizardo2022solving,herrera2024optimal,bayraktar2025deep}). However, these contributions are restricted to \textit{single} optimal stopping problems.

\noindent\textbf{The Multiple Stopping Challenge.}
In this paper, we address the numerical resolution of the \textit{multiple} optimal stopping problem, where one seeks to stop multiple stochastic processes by assigning a distinct stopping time to each. Applications of this framework are manifold, including the valuation of swing options in energy markets, the management of hydroelectric power plants, optimal liquidation strategies with multiple assets, and real options with multiple exercise rights.
Such problems were introduced in a general continuous-time framework by \citeauthor*{kobylanski2011optimal} \cite{kobylanski2011optimal}, proving that they reduce to a recursive sequence of single stopping problems with random horizons; the framework was recently extended to nonlinear expectations in \cite{grigorova2025non}. Related theoretical studies include asymptotic versions~\cite{talbi2023dynamic, talbi2023viscosity,talbi2024finite} and specific applications like swing options~\cite{carmona2008optimal}.
The numerical literature for multiple stopping is far more limited. While classical methods exist (e.g., \cite{meinshausen2004monte,schoenmakers2012pure}), \emph{deep learning} approaches are scarce. To the best of our knowledge, the primary existing deep learning method is that of \citeauthor*{han2023new} \cite{han2023new}, who extended the logic of \cite{becker2019deep} to multiple stopping for option pricing. On a related note, \cite{magninolearning} applied deep learning to mean-field optimal stopping by reformulating it as a mean-field Markov decision process. However the theoretical foundations remain very limited.

\noindent\textbf{Contributions.}
In this work, we propose a rigorous deep learning framework for solving multiple optimal stopping problems in high dimensions. Our main contributions are summarized as follows:

\begin{itemize}
    \item \textbf{Efficient Algorithm:} We develop a backward induction algorithm where the value function is approximated by a single neural network. We introduce an auxiliary state variable encoding the status (active or stopped) of each component, which allows us to avoid training separate networks for different subsets of components, streamlining the learning process. Furthermore, we stop at most one component at each time step instead of doing an exhaustive search (Algorithm~\ref{alg:DBMSP}).
    \item \textbf{Theoretical Convergence Analysis:} We provide a rigorous convergence proof for our algorithm, taking into account the error induced by stopping at most one component per time step, the statistical sampling error, and the neural network approximation error (Theorems~\ref{thm:convergence} and~\ref{thm:cv-optimal-strategy}). In the case of diffusion processes, we further derive explicit error bounds that account for the Euler discretization (Proposition~\ref{prop:discrete-continuous2}).
    \item \textbf{Scalable Implementation:} We demonstrate the efficacy of our method on high-dimensional multiple American put options and nonlinear utility maximization problems (Section~\ref{sec:numerics}).
\end{itemize}

The paper is organized as follows. Section~\ref{sec:setting} defines the problem and derives the dynamic programming principle. Section~\ref{sec:main} details the algorithm and establishes the convergence results. In Section~\ref{sec:diffusion}, we analyze the specific case of diffusion processes (Euler scheme) and present numerical experiments. Technical proofs are gathered in Appendix~\ref{sec:appendix}.
\\

\noindent\textbf{Notation.} For any $n \in \mathbb{N}$, we write $[n] := \{0, 1, \dots, n\}$ and $[n]^* := [n]\setminus \{0\}.$

\section{The general discrete-time setting}\label{sec:setting}

We consider a discrete-time interacting particle system where the dynamics of each particle are influenced by the states of the other particles until it is stopped. Let $p \in \mathbb{N}^*$ and $(\Omega, \mathcal{F}, \mathbb{P})$ be a probability space endowed with a filtration $\mathbb{F} := \{\mathcal{F}_n\}_{n \in [p]}$, with $\mathcal{F}_p = \mathcal{F}$. We denote by $\mathcal{T}_p$ the set of $[p]$-valued $\mathbb{F}$-stopping times, and by $\mathcal{T}_p^N$ the set of $N$-tuples of elements of $\mathcal{T}_p$. \\

Let $\{ \epsilon_n \}_{n \in [p]}$ be a sequence of i.i.d.\ random variables such that $\epsilon_{n+1}$ is independent of $\mathcal{F}_n$ for all $n \in [p-1]$. Given $\boldsymbol{\tau} := (\tau^1, \dots, \tau^N) \in \mathcal{T}_p^N$, we denote by $\bI := (I^1, \dots, I^N)$ the corresponding vector of survival processes, defined by $I_n^k := \mathbb{1}_{n \le \tau^k}$ for all $n \in [p]$ and $k \in [N]^*$. We then consider the dynamics:
\begin{equation}\label{dynamics}
 \bX_{n+1} := \bX_n + F_{n}(\bX_n, \epsilon_{n+1})\bI_{n+1}, \quad \bX_0 \in \mathbb{R}^N,
\end{equation}
where $F : \mathbb{R}^N \times \mathbb{R}^d \longrightarrow \mathbb{R}^N$, and the product above is understood as the Hadamard (element-wise) product.

\subsection{The original problem} 

Let $c : [p] \times \mathbb{R}^N \times \{0,1\}^N \times \{0,1\}^N \to \mathbb{R}_+$ and $g : \mathbb{R}^{N} \to \mathbb{R}$ denote respectively the reward and terminal reward functions, such that:  
$$ 
    c_n(\bx, \bi, \bi) = 0 \quad \text{for all} \ (n, \bx, \bi) \in [p] \times \mathbb{R}^N \times \{0,1\}^N.
$$
The multiple optimal stopping problem consists of solving the following maximization problem:
\begin{equation}\label{original-pb}
 V_0 := \sup_{\boldsymbol{\tau} \in \mathcal{T}_p^N} \mathbb{E}\bigg[\sum_{n=0}^{p-1} c_n(\bX_n, \bI_n, \bI_{n+1}) + g(X_{\tau^1}^1, \dots, X_{\tau^N}^N)\bigg] = \sup_{\boldsymbol{\tau} \in \mathcal{T}_p^N} \mathbb{E}\bigg[\sum_{n=0}^{p-1} c_n(\bX_n, \bI_n, \bI_{n+1}) + g(\bX_p)\bigg], 
\end{equation}
where the second equality follows from the fact that $\bX$ is the vector of stopped processes, i.e., $X_{\tau^k}^k = X_p^k$ for all $k \in [N]^*$. This class of problems may be analyzed using the dynamic programming approach. To this end, we introduce the dynamic value function: 
\begin{equation}\label{dynamic-original-pb}
V_n(\bx, \bi) := \sup_{\boldsymbol{\tau} \in \mathcal{T}_{n, p}^N} \mathbb{E}\bigg[\sum_{k=n}^{p-1} c_k(\bX_k, \bI_k, \bI_{k+1}) + g(\bX_p) \,\Big|\, \bX_n = \bx, \bI_n = \bi\bigg],   
\end{equation}
where $\mathcal{T}_{n, p}^N$ denotes the set of $\{n, \dots, p\}$-valued $N$-tuples of $\mathbb{F}$-stopping times. Similarly to \citeauthor*{kobylanski2011optimal} \cite{kobylanski2011optimal} and \citeauthor*{talbi2024finite} \cite{talbi2024finite}, we have the following result:
\begin{proposition}\label{prop:DPP}
Assume the functions $c$ and $g$ are continuous. Then: \\
${\rm (i)}$ For all $(\bx, \bi) \in \mathbb{R}^N \times \{0,1\}^N$, we have:
\begin{align}\label{DPP}
    V_n(\bx, \bi) &= \max_{\{\bi ' \in \{0,1\}^N \ \text{s.t.} \ \bi ' \le \bi\}} \mathbb{E}\Big[ c_n(\bX_n, \bI_n, \bI_{n+1}) + V_{n+1}(\bX_{n+1}, \bI_{n+1}) \,\Big|\, \bX_n = \bx, \bI_n = \bi  \Big] \nonumber \\
    &= \max_{\{\bi ' \in \{0,1\}^N \ \text{s.t.} \ \bi' \le \bi\}} \mathbb{E}\Big[ c_n(\bx, \bi, \bi') + V_{n+1}(\bx + F_{n}(\bx, \epsilon_{n+1}) \bi', \bi' )  \Big], \nonumber
\end{align}
where the inequality $\bi' \le \bi$ is understood coordinate-wise.\\
${\rm (ii)}$ There exists an optimal stopping strategy $\boldsymbol{\tau}^* \in \mathcal{T}_{p}^N$ for~\eqref{original-pb}. 
\end{proposition}
\proof 
To prove ${\rm (i)}$, let us introduce the set $\mathbb{I}^N$ of $\mathbb{F}$-predictable, decreasing, $\{0,1\}^N$-valued processes such that $\bI_0 = \mathbf{1}$. It is clear that $\mathcal{T}_p^N$ and $\mathbb{I}^N$ are in bijection. Indeed, to each $\boldsymbol{\tau} = (\tau_1, \dots, \tau_N)$, we can associate $\bI = (I^1, \dots, I^N)$ by setting $I_n^k := \mathbb{1}_{n \le \tau_k} = 1 - \mathbb{1}_{n-1 \ge \tau_k}$, for $k \in [N]^*$. Since $\tau_k$ is a stopping time, it is clear that $I^k$, and therefore $\bI$, is a predictable process. Conversely, given $\bI \in \mathbb{I}^N$, we define the corresponding stopping times by setting $\tau_k := \min\{ n \ge 0 : I_{n+1}^k = 0 \} \wedge p$. \\
Now, let $\mathcal{A}^N$ be the set of $\mathbb{F}$-adapted processes taking values in $\{0,1\}^N$. Then $\mathbb{I}^N$ is in bijection with a subset of $\mathcal{A}^N$ representing binary controls. Indeed, given $\bI_n = (\mathbb{1}_{n \le \tau_1}, \dots, \mathbb{1}_{n \le \tau_N})$, we have, by setting $\alpha_n^k := \mathbb{1}_{\{\tau_k = n\}}$ for $k \in [N]^*$, that $I_n^k = \prod_{j = 1}^{n-1} \alpha_j^k.$
Therefore, the multiple stopping problem~\eqref{dynamic-original-pb} may be rewritten as a standard control problem on $\mathcal{A}$. Since $g$ is continuous, the dynamic programming principle for this problem writes:
$$ V_n(\bx, \bi) = \max_{\mathbf{a} \in \{0,1\}^N} \mathbb{E}\Big[c_n(\bX_n, \bI_n, \bI_{n+1}) + V_{n+1}(\bX_{n+1}, \bI_{n+1}) \,\Big|\, \bX_n = \bx, \bI_n = \bi\Big], $$
where $\bX_{n+1} = \bx + F_n(\bx, \epsilon_{n+1})\bI_{n+1}$ and $\bI_{n+1} = \bi \mathbf{a}$. The desired result is obtained by observing that:
$$ \{ \bi \mathbf{a} : \mathbf{a} \in \{0,1\}^N\} = \{ \bi' \in \{0,1\}^N : \bi' \le \bi \}, $$
and that this set is finite; therefore, the supremum is a maximum.  \\
To prove ${\rm{(ii)}}$, we recursively construct the process $\bI^*$ defined by 
$$
    \bI_n^* = \bi \quad \text{and} \quad \bI_{n+1}^* \in {\rm{argmax}}_{\bi' \le \bI_n^*} \mathbb{E}\big[c_n(\bx, \bi, \bi') + V_{n+1}(\bX_{n+1}, \bi')\big].
$$
Then, $\boldsymbol{\tau}^*$ is the $N$-tuple of stopping times associated with the $N$-tuple of survival processes $\bI^*$. 
\qed 

Intuitively, Proposition~\ref{prop:DPP} implies that the multiple optimal stopping problem can be reduced to a recursive sequence of standard stopping problems: at each time $n$, one decides which agents will be stopped. This decision is encapsulated in the vector $\bi'$: if $i_k' = 0$, then the $k$-th agent is stopped; otherwise, they continue. Note that at every time $n$, given a state vector $\bi \in \{0,1\}^N$, one must examine all combinations of $\bi' \le \bi$ to decide which particles should be stopped. As will be seen in Section~\ref{sec:main}, this may induce a computational cost with exponential growth in $N$. We therefore examine an alternative multiple stopping problem which considerably reduces this cost. 

\subsection{The alternative problem} 

Introduce the following set of $N$-tuples of stopping times, in which no two elements of the tuple can be equal unless they are equal to the terminal time:
$$ \tilde{\mathcal{T}}_p^N := \Big\{ \boldsymbol{\tau} = (\tau_1, \dots, \tau_N) \in \mathcal{T}_p^N : \tau_k = \tau_l \implies \tau_k = p \quad \text{for all} \ k \neq l \Big\}. $$
We then define the alternative multiple optimal stopping problem: 
\begin{equation}
\label{alternative-pb}
     \tilde{V}_0 := \max_{\boldsymbol{\tau} \in \tilde{\mathcal{T}}_p^N} \mathbb{E}\bigg[\sum_{n=0}^{p-1} c_n(\bX_n, \bI_n, \bI_{n+1}) + g(X_{\tau^1}^1, \dots, X_{\tau^N}^N)\bigg] = \max_{\boldsymbol{\tau} \in \tilde{\mathcal{T}}_p^N} \mathbb{E}\bigg[\sum_{n=0}^{p-1} c_n(\bX_n, \bI_n, \bI_{n+1}) + g(\bX_p)\bigg].
\end{equation}
In this problem, at each time $n \in [p-1]$, one can stop (at most) one particle. Similarly to the original problem, we introduce a dynamic version of the value function in~\eqref{alternative-pb}:
 \begin{equation}\label{dynamic-alternative-pb}
     \tilde{V}_n(\bx, \bi) := \max_{\tilde{\boldsymbol{\tau}} \in \tilde{\mathcal{T}}_{n, p}^N} \mathbb{E}\bigg[\sum_{k=n}^{p-1} c_k(\bX_k, \bI_k, \bI_{k+1}) + g(\bX_p) \,\Big|\, \bX_n = \bx, \bI_n = \bi \bigg], 
 \end{equation}
where $\tilde{\mathcal{T}}_{n, p}^N$ denotes the set of $\{n, \dots, p\}$-valued $N$-tuples of $\tilde{\mathcal{T}}_p^N$. We then have the following dynamic programming principle:
\begin{proposition}\label{prop:DPP2} 
    Assume $g$ is continuous. Then:\\
    ${\rm{(i)}}$ For every $(\bx, \bi) \in \mathbb{R}^N \times \{0,1\}^N$, we have: 
    \begin{equation}
        \tilde{V}_n(\bx, \bi) = \max_{\ell \in [N]} \mathbb{E}\Big[ c_n(\bx, \bi, \bi^{-\ell}) +  \tilde{V}_{n+1}(\bx + F_{n}(\bx, \epsilon_{n+1}) \bi^{-\ell}, \bi^{-\ell} )  \Big], \nonumber
    \end{equation}
    where $\bi^{-\ell} := \bi \cdot (1 - \mathbf{e}_{\ell})$ for $\ell \in [N]^*$, and $\bi^{-0} := \bi$, with $(\mathbf{e}_1, \dots, \mathbf{e}_N)$ denoting the canonical basis of $\mathbb{R}^N$. \\
    ${\rm{(ii)}}$ There exists an optimal stopping strategy $\tilde{\boldsymbol{\tau}}^* \in \tilde{\mathcal{T}}_{n,p}^N$ for the problem~\eqref{dynamic-alternative-pb}. 
\end{proposition}
\proof 
The proof follows the same path as the proof of Proposition~\ref{prop:DPP}, replacing $\mathcal{A}^N$ with the set $\widetilde{\mathcal{A}}^N$ of $\mathbb{F}$-adapted processes taking their values in $\{ \mathbf{1} - \mathbf{e}_\ell : \ell \in [N] \}$.
\qed \\

According to Proposition~\ref{prop:DPP2}, at each step $n$, we must choose which particle to stop (if any). This reduces to choosing an index $\ell \in [N]$, with $\ell = 0$ representing the case where no particle is stopped. Compared with the original problem, we trade a possible exponential cost in $N$ for a linear cost in $N$. However, this cost reduction entails an additional error due to the fact that in the new problem, we cannot stop several particles simultaneously. In Section~\ref{sec:diffusion}, we analyze this error in the context of discretized diffusion processes.

\section{Main results}\label{sec:main}

\subsection{The original algorithm}

The first algorithm is based on Proposition~\ref{prop:DPP} and directly approximates the original problem~\eqref{original-pb}. The central idea is as follows: assuming the function $V_{n+1}$ has been appropriately approximated at time $n+1$, we compute the function $V_n$ in two steps: 
\begin{enumerate}
\item First, we approximate the function $U_n : (\bx, \bi) \mapsto \mathbb{E}\big[ V_{n+1}( \bx + F_n(\bx, \epsilon_{n+1})\bi, \bi)\big]$. To do this, we approximate $U_n$ using a dense neural network and employ the classical least squares characterization of conditional expectation combined with a Monte Carlo approach. This requires $M$ simulations $\{ (\bX_n^{(m)}, \bI_n^{(m)}, \epsilon_{n+1}^{(m)})\}_{1 \le m \le M }$ drawn from a distribution $\nu = \mu_X \otimes \mu_I \otimes \mu_\epsilon$, which does not depend on the current time $n$. 
\item Then, given $\bi \in \{0,1\}^N$, the function $V_n(\cdot, \bi)$ is defined as: 
$$V_n(\cdot, \bi) = \max_{\bi' \le \bi} \big[ c_n(\bx, \bi, \bi') + U_n(\bx, \bi') \big].$$ 
\end{enumerate}

In what follows, we denote by $\xi_M := \{ (\bX_n^{(m)}, \bI_n^{(m)}, \epsilon_{n+1}^{(m)}) : 1 \le m \le M, 0 \le n \le p - 1\}$ the full set of simulations, and by $\hat U_n^{\xi_M}$ and $\hat V_n^{\xi_M}$ the neural network approximations of the functions $U_n$ and $V_n$, respectively. The method is summarized in Algorithm~\ref{alg:DBMSE} and referred to as Exhaustive Deep Backward Multiple Stopping (E-DBMS for short). Notice that $\{\hat{V}^{\xi_M}_n\}_{n=0}^{p-1}$ and policies $\{\bI^{\xi_M}_{n+1}\}_{n=0}^{p-1}$ are based on the neural networks $\{\hat{U}^{\xi_M}_n\}_{n=0}^{p-1}$ so in practice, we can return the latter and compute the former two quantities. If needed, one can also train auxiliary neural networks for $\{\hat{V}^{\xi_M}_n\}_{n=0}^{p-1}$ and policies $\{\bI^{\xi_M}_{n+1}\}_{n=0}^{p-1}$, which increases the training time but reduces the inference time.

\begin{algorithm}
\caption{Deep Backward Multiple Stopping with Exhaustive search (E-DBMS)}\label{alg:DBMSE}
\begin{algorithmic}[1]
\State {\bf Input:} Training data $\xi_M = \{ (\bX_n^{(m)}, \bI_n^{(m)}, \epsilon_{n+1}^{(m)}) : 1 \le m \le M, 0 \le n \le p - 1\}$, hypothesis space $\mathcal{V}$ (neural networks).
\State {\bf Output:} Approximated value functions $\{\hat{V}_n\}_{n=0}^{p-1}$ and policies $\{\bI_{n+1}\}_{n=0}^{p-1}$.
\State Initialize $\hat V_p^{\xi_M}(\bx, \bi) \gets g(\bx)$ for all $\bx, \bi$.
\For{$n = p-1, \dots, 0$}
    \State \textbf{Step 1: Regression.} Compute $\hat U_n^{\xi_M}$ by minimizing the empirical squared loss:
    $$ 
        \hat U_n^{\xi_M} \gets \underset{\phi \in \mathcal{V}}{\mathrm{argmin}} \frac{1}{M} \sum_{m=1}^M \bigg| \phi(\bX_n^{(m)}, \bI_n^{(m)}) - \hat V_{n+1}^{\xi_M}\Big(\bX_n^{(m)} + F_n(\bX_n^{(m)}, \epsilon_{n+1}^{(m)})\bI_n^{(m)}, \bI_n^{(m)}\Big)  \bigg|^2. 
    $$
    \State \textbf{Step 2: Maximization.} Define $\hat V_n^{\xi_M}$ as the increasing envelope of $\hat U_n$ with respect to $\bi$ {}{} via the dynamic programming principle:
    $$ \hat V_n^{\xi_M}(\bx, \bi) \gets \max_{\bi' \in \{0,1\}^N : \bi' \le \bi}  \Big[ c_n(\bx, \bi, \bi') + \hat U_n^{\xi_M}(\bx, \bi') \Big]. $$
    \State Define the optimal policy update map:
    $$ 
        \bI_{n+1}^{\xi_M}(\bx, \bi) \in \underset{\bi' \in \{0,1\}^N : \bi' \le \bi}{\mathrm{argmax}}  \Big[ c_n(\bx, \bi, \bi') + \hat U_n^{\xi_M}(\bx, \bi') \Big]. 
    $$
\EndFor
\State \Return $\{\hat{V}^{\xi_M}_n\}_{n=0}^{p-1}$ and policies $\{\bI^{\xi_M}_{n+1}\}_{n=0}^{p-1}$
\end{algorithmic}
\end{algorithm}

An important drawback of Algorithm~\ref{alg:DBMSE} is that the maximization step is done through exhaustive search, which is computationally expensive, as we have:
$$\mathrm{Card}\big(\{\bi' \in \{0,1\}^N : \bi' \le \bi \}\big) = 2^{| \bi |_1}.$$ 
Computing the maximum over this set therefore implies a computational cost of order $\mathcal{O}(2^N)$. This motivates the introduction of an alternative algorithm based on the alternative multiple stopping problem~\eqref{alternative-pb}.

\subsection{The alternative algorithm}

The second algorithm is based on Proposition~\ref{prop:DPP2}. It directly approximates the alternative problem~\eqref{alternative-pb} and, consequently, approximates the original problem~\eqref{original-pb} in the context of discretized stochastic differential equations (see Proposition~\ref{prop:error_V_tildeV}). The procedure is as follows: 
\begin{enumerate}
\item First, we approximate the function $U_n : (\bx, \bi) \mapsto \mathbb{E}\big[ V_{n+1}( \bx + F_n(\bx, \epsilon_{n+1})\bi, \bi)\big]$ similarly to the previous algorithm. 
\item Then, given $\bi \in \{0,1\}^N$, the function $V_n(\cdot, \bi)$ is defined as: 
$$V_n(\cdot, \bi) = \max_{\ell \in [N]} \big[ c_n(\bx, \bi, \bi^{-\ell}) + U_n(\bx, \bi^{-\ell}) \big].$$
\end{enumerate}

The method is summarized in Algorithm~\ref{alg:DBMSP} and referred to as partial DBMS (P-DMBS).

\begin{algorithm}
\caption{Deep Backward Multiple Stopping with Partial stopping (P-DBMS)}\label{alg:DBMSP}
\begin{algorithmic}[1]
\State {\bf Input:} Training data $\xi_M$, hypothesis space $\mathcal{V}$ (neural networks).
\State {\bf Output:} Approximated value functions $\{\hat{V}_n\}_{n=0}^{p-1}$ and policies $\{\bI_{n+1}\}_{n=0}^{p-1}$.
\State Initialize $\hat V_p^{\xi_M}(\bx, \bi) \gets g(\bx)$ for all $(\bx, \bi)$.
\For{$n = p-1, \dots, 0$}
    \State \textbf{Step 1: Regression.} Compute $\hat U_n^{\xi_M}$ by minimizing the empirical squared loss:
    $$ \hat U_n^{\xi_M} \gets \underset{\phi \in \mathcal{V}}{\mathrm{argmin}} \frac{1}{M} \sum_{m=1}^M \bigg| \phi(\bX_n^{(m)}, \bI_n^{(m)}) - \hat V_{n+1}^{\xi_M}\Big(\bX_n^{(m)} + F_n(\bX_n^{(m)}, \epsilon_{n+1}^{(m)})\bI_n^{(m)}, \bI_n^{(m)}\Big)  \bigg|^2. $$
    \State \textbf{Step 2: Maximization.} Define $\hat V_n^{\xi_M}$ by maximizing over single stopping decisions:
    $$ \hat V_n^{\xi_M}(\bx, \bi) \gets \max_{\ell \in [N]} \Big[ c_n(\bx, \bi, \bi^{-\ell}) + \hat U_n^{\xi_M}(\bx, \bi^{-\ell}) \Big]. $$
    \State Define the optimal policy update map: 
    $$ \bI_{n+1}^{\xi_M}(\bx, \bi) \gets \bi^{-\ell^{\xi_M}(\bx, \bi)}, \quad \text{where } \ell^{\xi_M}(\bx, \bi) \in \underset{\ell \in [N]}{\mathrm{argmax}} \Big[ c_n(\bx, \bi, \bi^{-\ell}) + \hat U_n^{\xi_M}(\bx, \bi^{-\ell}) \Big]. $$
\EndFor
\State \Return $\{\hat{V}^{\xi_M}_n\}_{n=0}^{p-1}$ and policies $\{\bI^{\xi_M}_{n+1}\}_{n=0}^{p-1}$
\end{algorithmic}
\end{algorithm}

\subsection{The convergence results}

In order to analyze the convergence of the algorithms, we shall restrict the neural networks to the following class of functions:
$$ 
    \cN_M := \Big\{ f : \dbR^N \times \dbR^{K_M(2+N)+1} \ni (\bx, \th) \mapsto \sum_{j = 1}^{K_M} \a_j \si(\b_j \cdot \bx + \g_j) + \a_0, \ \mbox{with} \ \th := (\a_j, \b_j, \g_j)_{j} \Big\},
$$ 
where $\si : \dbR \to \dbR$ is some activation function, and $\{K_M\}_{M \ge 0}$ is such that 
$$
    \d_M := \frac{K_M}{M} \to 0 \ \mbox{as} \ M \to \infty. 
$$ 
We shall also use the following notation: given two sequences of variables $Y_M$ and $Z_M$, $M \ge 0$, we say that $Y_M = \cO_\dbP(Z_M)$ if there exists a constant $C \ge 0$ such that $\dbP\big( | Y_M | \le C | Z_M |\big) \to 1$ as $M \to \infty$. 

Since we essentially use neural networks to approximate conditional expectations, we also recall the following result from \citeauthor{kohler2006nonparametric} \cite[Corollary 1]{kohler2006nonparametric}, which will be useful to analyze the convergence of the algorithm in the sense of $\cO_\dbP$:
\begin{lemma}\label{lem:conditional-expectation}
Let $(X_i, Y_i)_{1 \le i \le M}$ be a sequence of i.i.d.\ $\dbR^d \times \dbR$-valued random variables. Introduce the measurable functions:
$$ 
    \m(x) := \dbE[Y_1 | X_1 = x], \quad \m_M \in \underset{\phi \in \cN_M}{\mathrm{arg min}} \sum_{i=1}^M | \phi(X_i) - Y_i |^2. 
$$
We assume that there exist two positive constants $\sigma, \lambda$ such that:
$$ 
    \dbE\Big[ \exp\Big(\frac{(Y_1-m(X_1))^2}{\si^2}\Big) \Big| X_1\Big] \le \l. 
$$
Then we have:
$$ 
    \dbE\Big[ \Big|\m_M(X_1) - \m(X_1) \Big|^2  \Big] = \cO_\dbP\Big( \d_M + \inf_{\phi \in \cN_M}  \dbE\Big[ \Big| \phi(X_1) - m(X_1) \Big|^2 \Big]  \Big). 
$$
\end{lemma}

Before stating the main result, we introduce the metric used to measure the error: for $f : \dbR^N \times \{0,1\}^N \to \dbR$, we denote:
$$ \lVert f \rVert_{2,\infty}^{\xi_M} := \dbE\Big[ \max_{\bi \in \{0,1\}^N} \big| f(\bX_0, \bi) \big|^2 \Big| \xi_M \Big]^{1/2},  $$
where $\dbP \circ \bX_0^{-1} = \m$.

\begin{assumption}\label{assum:convergence}
Assume that $\dbP \circ \bX_n^{-1} = \m$. Then, for all $\bi \in \{0,1\}^N$, the random variable $\bX_{n+1} = \bX_n + F_n(\bX_n, \e_{n+1})\bi$ admits a bounded density with respect to $\m$ conditionally on $\bX_n$, denoted $\bx \mapsto h_n(\bx ; \bX_n, \bi)$.
\end{assumption}

\begin{theorem}\label{thm:convergence}
Let Assumption~\ref{assum:convergence} hold. \\
${\rm (i)}$ Let $\hat V^{\xi_M}$ be the function resulting from Algorithm~\ref{alg:DBMSE}. Then we have, as $M \to \infty$:
\bea
    \lVert \hat V_0^{\xi_M} - V_0 \rVert_{2, \infty}^{\xi_M} = \cO_\dbP\Big(\d_M + \sup_{0 \le k \le p}\inf_{\phi \in \cN_M} \lvert \phi - U_k \rvert_{2, \infty}^{\xi_M} \Big), \nonumber
\eea
for some $\d_M \to 0$ as $M \to \infty$.

${\rm (ii)}$ Let $\hat V^{\xi_M}$ be the function resulting from Algorithm~\ref{alg:DBMSP}. Then we have, as $M \to \infty$:
\bea
    \lVert \hat V_0^{\xi_M} - \tilde V_0 \rVert_{2, \infty}^{\xi_M} = \cO_\dbP\Big(\d_M + \sup_{0 \le k \le p}\inf_{\phi \in \cN_M} \lvert \phi - U_k \rvert_{2, \infty}^{\xi_M} \Big), \nonumber
\eea
for some $\d_M \to 0$ as $M \to \infty$.
\end{theorem}
\proof We only write the proof of (i), as (ii) proceeds exactly from the same arguments. For $n \in [p]$, introduce the following function:
\bea 
\bar U_n^{\xi_M}(\bx, \bi) &:=& \dbE\big[\hat V_{n+1}^{\xi_M}\big(\bx + F_n(\bx, \e_{n+1})\bi, \bi \big)\big], \nonumber
\eea
for all $(\bx, \bi) \in \dbR^N \times \{0,1\}^N$, and we recall that $U_n$ denotes the function defined by:
\bea
U_n(\bx, \bi) &:=& \dbE\big[V_{n+1}\big(\bx + F_n(\bx, \e_{n+1})\bi, \bi \big)\big]. \nonumber 
\eea
Then we have, by definition of Algorithm~\ref{alg:DBMSE} and Proposition~\ref{prop:DPP}:
\begin{align*}
    \lVert \hat V_n^{\xi_M} - V_n \rVert_{2, \infty}^{\xi_M} 
    &\le \dbE\Big[ \max_{\bi \in \{0,1\}^N} \big| \max_{\bi' \le \bi} \hat U_n^{\xi_M}(X_n, \bi') - \max_{\bi' \le \bi}  U_n(X_n, \bi') \big|^2 \Big| \xi_M \Big]^{1/2} 
    \\
    &\le \lVert \hat U_n^{\xi_M} - U_n \rVert_{2, \infty}^{\xi_M} 
    \\
    &\le \lVert \hat U_n^{\xi_M} - \bar U_n^{\xi_M} \rVert_{2, \infty}^{\xi_M} + \lVert \bar U_n^{\xi_M} - U_n \rVert_{2, \infty}^{\xi_M}.
\end{align*}
Now, introduce the subset of $\Omega$:
$$ A_n^M := \Big\{ \lVert \hat U_n^{\xi_M} - \bar U_n^{\xi_M} \rVert_{2, \infty}^{\xi_M} \le C_n\Big(\d_M + \inf_{\phi \in \cN_M} \lVert \phi - \bar U_n^{\xi_M} \rVert_{2, \infty}^{\xi_M}\Big)\Big\}, $$
where $\d_M \to 0$ and $C_n$ is such that $\dbP(A_n^M) \to 1$ as $M \to \infty$, by Lemma~\ref{lem:conditional-expectation}. On this set, we have:
\bea\label{ineq1-convergence}
\lVert \hat V_n^{\xi_M} - V_n \rVert_{2, \infty}^{\xi_M} &\le& C_n\Big(\d_M + \inf_{\phi \in \cN_M} \lVert \phi - \bar U_n^{\xi_M} \rVert_{2, \infty}^{\xi_M}\Big) + \lVert \bar U_n^{\xi_M} - U_n \rVert_{2, \infty}^{\xi_M} \nonumber \\
&\le& C_n\Big(\d_M + \inf_{\phi \in \cN_M} \lVert \phi - U_n \rVert_{2, \infty}\Big) + (1+C_n)\lVert \bar U_n^{\xi_M} - U_n \rVert_{2, \infty}^{\xi_M}.
\eea
Now, observe that:
\begin{align*}
    \lVert \bar U_n^{\xi_M} - U_n \rVert_{2, \infty}^{\xi_M} 
    &\le \dbE\Big[\max_{\bi \in \{0,1\}^N} \big| \dbE\big[\hat V_{n+1}^{\xi_M}(\bX_{n+1}, \bi) | \bX_n, \xi_M \big] - \dbE\big[V_{n+1}(\bX_{n+1}, \bi) | \bX_n, \xi_M \big] \big|^2 \Big| \xi_M\Big]^{1/2} 
    \\
    &\le \dbE\Big[\dbE\Big[ \max_{\bi \in \{0,1\}^N}  \big| \hat V_{n+1}^{\xi_M}(\bX_{n+1}, \bi) - V_{n+1}(\bX_{n+1}, \bi) \big|^2 | \bX_n\Big]  \Big| \xi_M\Big]^{1/2}
\end{align*}
with $\bX_{n+1} = \bX_n + F_n(\bX_n, \e_{n+1})\bi$. Now, using Assumption~\ref{assum:convergence}, we obtain:
\begin{align*}
\lVert \bar U_n^{\xi_M} - U_n \rVert_{2, \infty}^{\xi_M} \le& \dbE\Big[ \int_{\dbR^N} \max_{\bi \in \{0,1\}^N}  \big| \hat V_{n+1}^{\xi_M}(\bx, \bi) - V_{n+1}(\bx, \bi) \big|^2 h_n(\bx; \bX_n, \bi)\m(d\bx) | \xi_M \Big]^{1/2} \\
\le& \lVert h \rVert_\infty \dbE\Big[ \int_{\dbR^N} \max_{\bi \in \{0,1\}^N}  \big| \hat V_{n+1}^{\xi_M}(\bx, \bi) - V_{n+1}(\bx, \bi) \big|^2 \m(d\bx) | \xi_M \Big]^{1/2} \\
=& \lVert h \rVert_\infty \lVert \hat V_{n+1}^{\xi_M} - V_{n+1} \rVert_{2,\infty}^{\xi_M}. 
\end{align*}
Plugging this into~\eqref{ineq1-convergence}, we have on $A_n^M$:
\begin{align*}
 \lVert \hat V_n^{\xi_M} - V_n \rVert_{2, \infty}^{\xi_M} \le& C_n\Big(\d_M + \inf_{\phi \in \cN_M} \lVert \phi - U_n \rVert_{2, \infty}\Big) + (1+C_n)\lVert h \rVert_\infty \lVert \hat V_{n+1}^{\xi_M} - V_{n+1} \rVert_{2,\infty}^{\xi_M} \\
 \le& C\Big(\d_M + \sup_{0 \le k \le p}\inf_{\phi \in \cN_M} \lVert \phi - U_k \rVert_{2, \infty}\Big) + (1+C)\lVert h \rVert_\infty \lVert \hat V_{n+1}^{\xi_M} - V_{n+1} \rVert_{2,\infty}^{\xi_M}, 
 \end{align*}
 with $C := \max_{n \in [p]} C_n$. Then, by induction, we have on $\cap_{n=0}^{p-1} A_n^M$, using the fact that $\hat V_p^{\xi_M} = g$:
 $$ \lVert \hat V_0^{\xi_M} - V_0 \rVert_{2, \infty}^{\xi_M} \le C(1+C)^p\lVert h \rVert_\infty^p\Big(\d_M + \sup_{0 \le k \le p}\inf_{\phi \in \cN_M} \lVert \phi - U_k \rVert_{2, \infty} \Big).   $$
We conclude the proof by observing that $\dbP\big[(\cap_{n=0}^{p-1} A_n^M)^c\big] \le \sum_{n=0}^p \dbP[(A_n^M)^c] \to 0$ as $M \to \infty$. 
\qed

\begin{theorem}\label{thm:cv-optimal-strategy}
{\rm (i)} Let $\bI^{\xi_M}$ be the stopping strategy provided by Algorithm~\ref{alg:DBMSE}, and denote for $n \in [p]$:
$$ 
    J_n^{\xi_M} := \dbE\Big[ \sum_{k=n}^{p-1} c_k(\bX_k, \bI_k^{\xi_M}, \bI_{k+1}^{\xi_M}) + g(\bX_p) \big],
$$
where the dynamics of $\bX$ is controlled by $\bI^M$. Then we have:
$$ \lVert J_0^{\xi_M} - V_0 \rVert_{2,\infty}^{\xi_M} = \cO_\dbP\Big(\d_M + \sup_{0 \le k \le p} \inf_{\phi \in \cN_M} \lVert \phi - U_k \rVert_{2, \infty}\Big),$$
with $\d_M \to 0$ as $M \to \infty$. \\
{\rm (ii)} Let $\bI^{\xi_M}$ the stopping strategy provided by Algorithm~\ref{alg:DBMSP}, and define $J_n^{\xi_M}$ as above. Then we have:
$$ \lVert J_0^{\xi_M} - \tilde V_0 \rVert_{2,\infty}^{\xi_M} = \cO_\dbP\Big(\d_M + \sup_{0 \le k \le p} \inf_{\phi \in \cN_M} \lVert \phi - U_k \rVert_{2, \infty}\Big),$$
with $\d_M \to 0$ as $M \to \infty$.
\end{theorem}
\proof
For simplicity, we write the proof for $c = 0$. We only detail the argument for (i), as (ii) is proved in the very same way. First, observe that:
\begin{equation}\label{1stIneq}
    \D_n := \lVert J_n^{\xi_M} - V_n \rVert_{2,\infty}^{\xi_M} \le \lVert \hat V_n^{\xi_M} - V_n \rVert_{2,\infty}^{\xi_M} + \lVert \hat V_n^{\xi_M} - J_n^{\xi_M} \rVert_{2,\infty}^{\xi_M}.
\end{equation}
Yet, we have:
\begin{align*}
    \lVert \hat V_n^{\xi_M} - J_n^{\xi_M} \rVert_{2,\infty}^{\xi_M} &= \dbE\Big[\max_{\bi \in \{0,1\}^N} \Big| \hat V_n^{\xi_M}(\bX_n, \bi) - J_n^{M}(\bX_n, \bi)  \Big|^2 | \xi_M \Big]^{1/2} \\
    =& \dbE\Big[\max_{\bi \in \{0,1\}^N} \Big| \hat U_n^{\xi_M}(\bX_n, \bI_{n+1}^{\xi_M}) - \dbE\Big[J_{n+1}^{\xi_M}(\bX_{n+1}, \bI_{n+1}^{\xi_M}) | \bX_n \Big] \Big|^2 | \xi_M \Big]^{1/2}  \\
    \le& \lVert \hat U_n^{\xi_M} - \bar U_n^{\xi_M} \rVert_{2,\infty}^{\xi_M} + \dbE\Big[\max_{\bi \in \{0,1\}^N} \Big| \bar U_n^{\xi_M}(\bX_n, \bI_{n+1}^{\xi_M}) - \dbE\Big[J_{n+1}^{\xi_M}(\bX_{n+1}, \bI_{n+1}^{\xi_M}) | \bX_n \Big] \Big|^2 | \xi_M \Big]^{1/2} \\
    \le& \lVert \hat U_n^{\xi_M} - \bar U_n^{\xi_M} \rVert_{2,\infty}^{\xi_M} + \dbE\Big[\max_{\bi \in \{0,1\}^N} \Big| \bar U_n^{\xi_M}(\bX_n, \bI_{n+1}^{\xi_M}) - \dbE\Big[V_{n+1}(\bX_{n+1}, \bI_{n+1}^{\xi_M}) | \bX_n \Big] \Big|^2 | \xi_M \Big]^{1/2} \\
    &+ \dbE\Big[\max_{\bi \in \{0,1\}^N} \Big| \dbE\Big[V_{n+1}(\bX_{n+1}, \bI_{n+1}^{\xi_M}) | \bX_n \Big] - \dbE\Big[J_{n+1}^{\xi_M}(\bX_{n+1}, \bI_{n+1}^{\xi_M}) | \bX_n \Big] \Big|^2 | \xi_M \Big]^{1/2} \\
    \le& \lVert \hat U_n^{\xi_M} - \bar U_n^{\xi_M} \rVert_{2,\infty}^{\xi_M} + \dbE\Big[\max_{\bi \in \{0,1\}^N} \dbE\Big[ \Big|  \hat V_{n+1}^{\xi_M}(\bX_{n+1}, \bI_{n+1}^{\xi_M}) - V_{n+1}(\bX_{n+1}, \bI_{n+1}^{\xi_M}) \Big|^2 | \bX_n \Big]   | \xi_M \Big]^{1/2} \\
    &+ \dbE\Big[\max_{\bi \in \{0,1\}^N}  \dbE\Big[ \Big| V_{n+1}(\bX_{n+1}, \bI_{n+1}^{\xi_M}) - J_{n+1}^{\xi_M}(\bX_{n+1}, \bI_{n+1}^{\xi_M}) \Big|^2 | \bX_n \Big]  | \xi_M \Big]^{1/2} \\
    \le& \lVert \hat U_n^{\xi_M} - \bar U_n^{\xi_M} \rVert_{2,\infty}^{\xi_M} + \lVert h \rVert_\infty \lVert \hat V_{n+1}^{\xi_M} - V_{n+1} \rVert_{2,\infty}^{\xi_M} + \lVert h \rVert_\infty \lVert \hat V_{n+1} - J_{n+1}^{\xi_M} \rVert_{2,\infty}^{\xi_M}, 
\end{align*}
where we used Assumption~\ref{assum:convergence} as in the proof of Theorem~\ref{thm:convergence}. Plugging this into~\eqref{1stIneq}, we obtain:
\begin{align*}
    \D_n \le \lVert \hat V_n^{\xi_M} - V_n \rVert_{2,\infty}^{\xi_M} + \lVert \hat U_n^{\xi_M} - \bar U_n^{\xi_M} \rVert_{2,\infty}^{\xi_M} + \lVert h \rVert_\infty \lVert \hat V_{n+1}^{\xi_M} - V_{n+1} \rVert_{2,\infty}^{\xi_M} + \lVert h \rVert_\infty \D_{n+1},
\end{align*}
from which we deduce the desired result from Theorem~\ref{thm:convergence}, Lemma~\ref{lem:conditional-expectation} and the fact that $\D_p = 0$.
\qed

\begin{remark}
It is worth noting that the error due to the neural network approximation in Theorems~\ref{thm:convergence} and~\ref{thm:cv-optimal-strategy} is expressed in terms of the auxiliary functions $\{U_n\}_{n \in [p]}$, rather than the value functions $\{V_n\}_{n \in [p]}$ as in \cite{hure2021deep}. This stems from the fact that, in our methodology, only the functions $U_n$ are directly approximated by neural networks, while the approximation of $V_n$ is obtained as a maximum of a family of neural networks, rather than being a neural network itself.
\end{remark}

\section{Application to diffusion processes}\label{sec:diffusion}

Our objective is to numerically compute the value function of the multiple optimal stopping problem for a $N$-dimensional stopped diffusion, as defined by \citeauthor*{talbi2024finite} \cite{talbi2024finite}. More precisely, we are interested in the following problem:
\bea\label{eq:multiple-OS}
V_0 &:=& \sup_{\boldsymbol{\tau} \in \cT_{[0,T]}^N} \mathbb{E}\Big[\sum_{0 \le s \le T} \bm{c}_s(\bX_s)\cdot (\bI_s - \bI_{s-}) + g\big(\bX_T \big)\Big] \\
&=& \sup_{\boldsymbol{\tau} \in \cT_{[0,T]}^N} \mathbb{E}\Big[\sum_{k=1}^N c^k_{\t_k}(\bX_{\t_k}) + g\big(X_{\tau_1}^1, \dots, X_{\tau_N}^N\big)\Big]. \nonumber
\eea
where $\cT_{[0,T]}^N$ denotes the set of $[0,T]$-valued $N$-tuples of stopping times and $\bX := (X^1, \dots, X^N)$ is the system of interacting stopped diffusions: 
\begin{equation}\label{continuous-time-dynamics}
 dX_t^k = I_t^k\left(b_k(t, \bX_t)dt + \sigma_k(t, \bX_t)dW_t^k + \sigma_0(t, \bX_t)dW_t^0\right) 
 \end{equation}
for all $k \in [N]^* := \{1, \dots, N\}$, with $I_t^k := \boldsymbol{1}_{\tau^k \ge t}$ and where the standard Brownian motions $W^0, \dots, W^N$ are independent. The following assumption will be in force throughout this Section:
\begin{assumption}\label{assum:coef-diffusion}
For $\phi \in \{\b_k, \si_k, \si_0, c^k, g\}, k \in [N]^*$, $\phi$, there exists a nonnegative constant $L \ge 0$ and $\beta \in (0,1]$ such that:
\begin{align*}
&| \phi(t, \bx) - \phi(t, \bx')| \le C|\bx - \bx'|, \\
&| \phi(t, \bx) - \phi(s, \bx)| \le C|t - s|^\b,
\end{align*}
for all $(t,s) \in [0,T]^2$ and $(\bx, \bx') \times \dbR^N \times \dbR^N$. 
\end{assumption}

\subsection{Discrete-time approximation}

Let $\pi := \{t_0, \dots, t_p\}$ be a partition of $[0,T]$, with $p \in \dbN^*$. For simplicity, we assume that all the subintervals have the same size, i.e., $t_k = k\frac{T}{p}$ for all $k \in [p]$. In what follows, we shall denote by $h := \frac{T}{p}$ the step of the partition $\pi$. In this paragraph, we denote by $\cT_p$ the set of $\pi$-valued stopping times, and by $\cT_p^N$ the set of $N$-tuples of $\pi$-valued stopping times. We then introduce the discrete time multiple optimal stopping problem:
\begin{equation}\label{eq:discrete-multiple-OS}
V_0^{h} := \sup_{\boldsymbol{\tau} \in \cT_p^N } \mathbb{E}\Big[\sum_{k=1}^N c^k_{\t_k}(\bX_{\t_k}^h) + g(X_{\tau_1}^{1,h}, \dots, X_{\tau_N}^{N,h})\Big],
\end{equation}
with $\bX^h := (X^{1,h}, \dots, X^{N,h})$ denotes the Euler scheme of $\bX$ for the partition $\pi$, that is:
\bea\label{eq:euler-scheme} 
    X_{t}^{k, h} = x_k + \int_0^t  I_{t^p(s)+h}^k b_k(t^p(s), \bX_{t^p(s)}^h)ds + \int_0^t  I_{t^p(s)+h}^k \si_k(t^p(s), \bX_{t^p(s)}^h)dW_s^k \\ + \int_0^t  I_{t^p(s)+h}^k \si_0(t^p(s), \bX_{t^p(s)}^h)dW_s^0 \nonumber
\eea
 for all $k \in [N]^*$ and $t \in [0,T]$, with $t_p(s)$ is the largest $t' \in \pi$ such that $t' \le s$. This dynamics correspond to the general dynamics~\eqref{dynamics} with: 
 $$ 
    F_n(\bx, \e_{n+1}) = \begin{pmatrix} b_1(n, \bx) \\ \dots \\ b_N(n, \bx) \end{pmatrix} h + \begin{pmatrix} \si_1(n, \bx)\e_{n+1}^{(1)} \\ \dots \\ \si_N(n, \bx)\e_{n+1}^{(N)} \end{pmatrix} + \si_0(n, \bx)\e_{n+1}^{(0)}\begin{pmatrix} 1 \\ \dots \\ 1 \end{pmatrix},$$ and $\e_{n+1}^{(k)} := W_{t_{n+1}}^k - W_{t_n}^k$ for all $k \in [N]$. 
 In the following Lemma (whose proof is relegated to the appendix), we estimate the error between $\bX$ and its Euler scheme, for stopping times taking their values in the discrete set $\{t_0, \dots, t_p\}$:
 \begin{lemma}\label{lem:Euler-error}
Let Assumption~\ref{assum:coef-diffusion} hold. For all $\bm{\t} \in \cT_p^N$, we have:
$$ 
    \sup_{\bm{\tau} \in \cT_p^N} \dbE\Big[ \sup_{t \in [0,T]} | \bX_t - \bX_t^h |^2 \Big] \le C_{T,N} h^{2\b \wedge 1}, 
$$
where the constant $C_{T,N}$ only depends on $T$ and $N$.
 \end{lemma}
This estimate implies the following estimate between the two value functions:
\begin{proposition}\label{prop:discrete-continuous}
Let Assumption~\ref{assum:coef-diffusion} hold. Then we have:
$$ |V_0 - V_0^h | \le C_{T,N}h^{\b \wedge 1/2}, $$
where the constant $C_{T,N}$ only depends on $T$ and $N$.
\end{proposition}
\proof
Introduce the value function of the optimal stopping problem of $\bX$ on the set of $\pi$-valued stopping times:
$$ V_0^p := \sup_{\bm{\tau} \in \cT_p^N} \dbE\Big[\sum_{k=1}^N c^k_{\t_k}(\bX_{\t_k}) + g\big(X_{\t_1}^1, \dots, X_{\t_N}^N)\Big]. $$
We have:
$$ |V_0 - V_0^h | \le |V_0 - V_0^p | + |V_0^p - V_0^h |.$$
Assumption~\ref{assum:coef-diffusion} and Lemma~\ref{lem:Euler-error} imply that:
$$ |V_0^p - V_0^h | \le Ch^{\b \wedge 1/2}. $$
for some constant $C \ge 0$. The inequality $V_0^p \le V_0$ is clear, as it simply comes from the fact that $\cT_p^N \subset \cT_{[0,T]}^N$. Now, for $\e >0$, let $\bm{\tau}^\e \in \cT_{[0,T]}^N$ be an $\e$-optimal policy for the problem~\eqref{eq:multiple-OS}.  We then define $\bar{\bm{\tau}}^\e = (\bar \tau_1^\e, \dots, \bar \tau_p^\e)$ as follows: for each $k \in [p]$, $\bar \tau_k^\e$ is the smallest $t_m \in \pi$ such that $\tau_k^\e \le t_m$. Since $\{ \bar \tau_k^\e \le t_m \} = \{ \tau_k^\e \in (t_{m-1}, t_m] \} \in \cF_{t_m}$, $\bar \tau_k^\e$ is a stopping time for the filtration $\dbF^p$, and therefore $\bar{\bm{\tau}}^\e \in \cT_p^N$. Then, by Assumption~\ref{assum:coef-diffusion} and the fact that $g$ is Lipschitz-continuous, Lemma~\ref{lem:continuity-tau} and the fact that $\bar{\bm{\tau}}^\e \in \cT_p^N$, we have:
\begin{align*}
    V_0 
    &\le \dbE\Big[\sum_{k=1}^N c^k_{\t_k^\e}(\bX_{\t_k^\e}) + g\big(X_{\t_1^\e}^1, \dots, X_{\t_N^\e}^N\big)\Big] + \e 
    \\
    &\le \dbE\Big[\sum_{k=1}^N c^k_{\bar \t_k^\e}(\bX_{\bar \t_k^\e}) + g\big(X_{\bar \t_1^\e}^1, \dots, X_{\bar \t_N^\e}^N\big)\Big]+ C_{T,N} h^{\b \wedge 1/2} + \e 
    \\
    &\le V_0^h + C_{T,N} h^{\b \wedge 1/2} + \e.
\end{align*}
We conclude by arbitrariness of $\e > 0$, and remarking that we typically have $h \le 1$.
\qed

\begin{remark}
If $g$ only depends on the empirical measure of $\bX$, one can show that the constant $C_{T, N}$ in fact does not depend on $N$.
\end{remark}

\subsection{Error for the alternative problem} 

In this paragraph, we consider the alternative problem
\begin{equation}\label{eq:discrete-multiple-OS2}
\tilde V_0^{h} := \sup_{\boldsymbol{\tau} \in \tilde \cT_p^N } \mathbb{E}\Big[\sum_{k=1}^N c^k_{\t_k}(\bX_{\t_k}^h) + g(X_{\tau_1}^{1,h}, \dots, X_{\tau_N}^{N,h})\Big],
\end{equation}
where agents may only be stopped one by one. As it is more challenging to compare $\tilde V_0^h$ to $V_0$ directly, we start by comparing $\tilde V_0^h$ to $V_0^h$:
\begin{proposition}\label{prop:error_V_tildeV}
Let Assumption~\ref{assum:coef-diffusion} hold. Then we have, for some constant $C_{T,N}$ depending on $T$ and $N$ only:
$$ 
    V_0^h - C_{T,N} h^{\beta \wedge 1/2}  \le \tilde{V}_0^h \le V_0^h.
$$
\end{proposition}
\proof 
The inequality $\tilde{V}_0^h \le V_0^h$ is clear, since $\tilde{\cT}_p^N \subset \cT_p^N$. For the other inequality, consider an optimal strategy $\boldsymbol{\tau}^* = (\tau_1^*, \dots, \tau_N^*)$ for~\eqref{eq:discrete-multiple-OS}, whose existence is guaranteed by Proposition~\ref{prop:DPP}, that is:
$$ 
    V_0^h = \dbE\Big[\sum_{k=1}^{N} c^k_{\t_k^*}(\bX_{\t_k^*}^{h,*}) + g(\bX_T^{h,*})\Big], 
$$
where $(\bX^{h,*}, \bI^*)$ denotes the process~\eqref{eq:euler-scheme} controlled by $\bm{\tau^*}$. 
By Lemma~\ref{lem:tildetau}, we may construct a strategy $\tilde{\boldsymbol{\tau}} = (\tilde \tau_1, \dots, \tilde \tau_N) \in \tilde{\cT}_p^N$ such that
$$ 
    | \tau_k^* - \tilde \tau_k | \le Nh, \ \mbox{a.s.}, 
$$ 
and we denote by $\tilde{\bX}^h$ the Euler scheme controlled by $\tilde{\boldsymbol{\tau}}$. By Assumption~\ref{assum:coef-diffusion}, we have the estimates:
\begin{align*}
    V_0^h \le&  \dbE\Big[ L \Big(\sum_{k=1}^N  (1 + | \bX_{\tau_k^{*}}^{h,*}|) |\tau_k^* - \tilde{\tau}_k |^\beta + |\bX_{\tau_k^{*}}^{h,*} - \tilde{\bX}_{\tilde{\tau}_k}^h |  + | \bX^{h,*}_T - \tilde{\bX}_T^h | \Big) + \sum_{k=1}^N c^k_{\tilde \t_k}(\tilde{\bX}_{\tilde{\t}_k}^h) +  g(\tilde{\bX}_T^h)  \Big] \\
    \le& L\Big( \big(1 + \sup_{n \in [p]} \dbE\big[ | \bX_{t_n}^{h,*} |]\big) h^\beta + C'_{T,N}\sqrt{\dbE\big[ | \bm{\tau} - \tilde{\bm{\tau}} | \big]} \Big) + \tilde{V}_0^h \\
    \le&  C_{T,N}h^{\b \wedge 1/2} + \tilde{V}_0^h,
\end{align*}
where we used Lemma~\ref{lem:error} for the second inequality.
\qed

Combining Propositions~\ref{prop:discrete-continuous} and~\ref{prop:error_V_tildeV}, we immediately deduce the following result:
\begin{proposition}\label{prop:discrete-continuous2}
Let Assumption~\ref{assum:coef-diffusion} hold. Then we have:
$$ |V_0 - \tilde{V}_0^h | \le C_{T,N}h^{\b \wedge 1/2}, $$
where the constant $C_{T,N}$ only depends on $T$ and $N$.
\end{proposition}

\subsection{Numerical implementation}
\label{sec:numerics}

We now present two examples to illustrate the performance of our approach.

\paragraph*{Multiple American Put.} We first test our algorithms on the following problem, which corresponds to the price of a basket of many American Puts, defined for all $\bx = (x_1, \dots, x_N) \in \dbR_+^N$ by:
$$ 
    V_0(\bx) := \sup_{\boldsymbol{\t} \in \cT_{[0,T]}^N} \dbE\Big[ \sum_{j=1}^N \big( K_j - S_{\t_j}^j \big)_+ \Big], 
$$
where $K_j \ge 0$, $S_t^j = x_j\exp\big([\m_j - \si_j^2/2]t + \si_j W_t^j)$, $j \in [N]^*$, with $W^1, \dots, W^N$ independent Brownian motions. Denoting by $S^{j,h}$ the Euler scheme of $S^j$, we also introduce:
$$ 
    V_0^h(\bx) := \sup_{\boldsymbol{\t} \in \cT_{p}^N} \dbE\Big[ \sum_{j=1}^N \big( K_j - S_{\t_j}^{j,h} \big)_+ \Big]. 
$$
Our objective is to compute $V_0^h$ through our deep learning algorithm and to compare to $V_0$, which can be easily approximated by finite difference once we observe that $V_0(\bx) = \sum_{j=1}^N V_0^j(x_j)$, with:
$$ 
    V_0^j(x_j) := \sup_{\t_j \in \cT_{[0,T]}} \dbE\Big[ \big( K_j - S_{\tau_j}^j \big)_+ \Big], 
$$
which corresponds to a single-agent optimal stopping problem.

\begin{figure}
\centering
\begin{subfigure}{.48\textwidth}
  \centering
    \includegraphics[width=.95\linewidth]{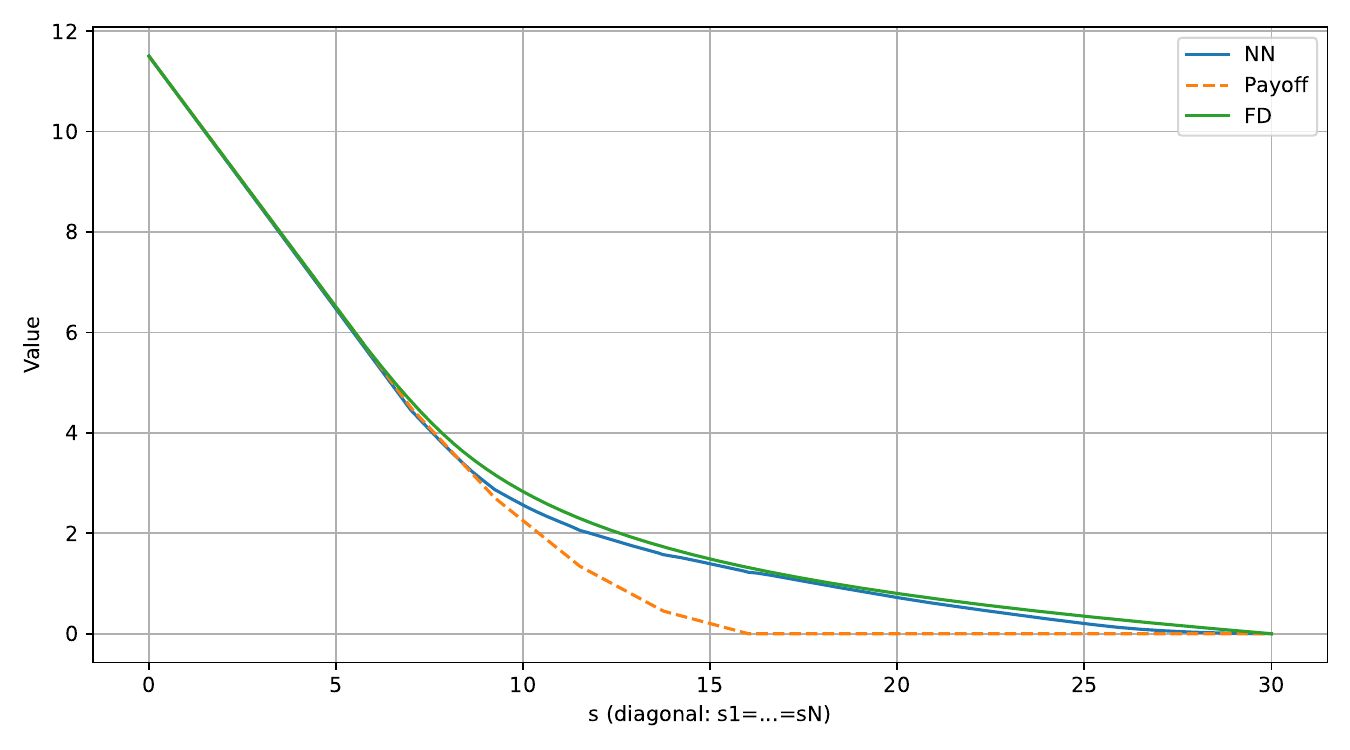}
\end{subfigure}
\centering
\begin{subfigure}{.48\textwidth}
  \centering
    \includegraphics[width=.95\linewidth]{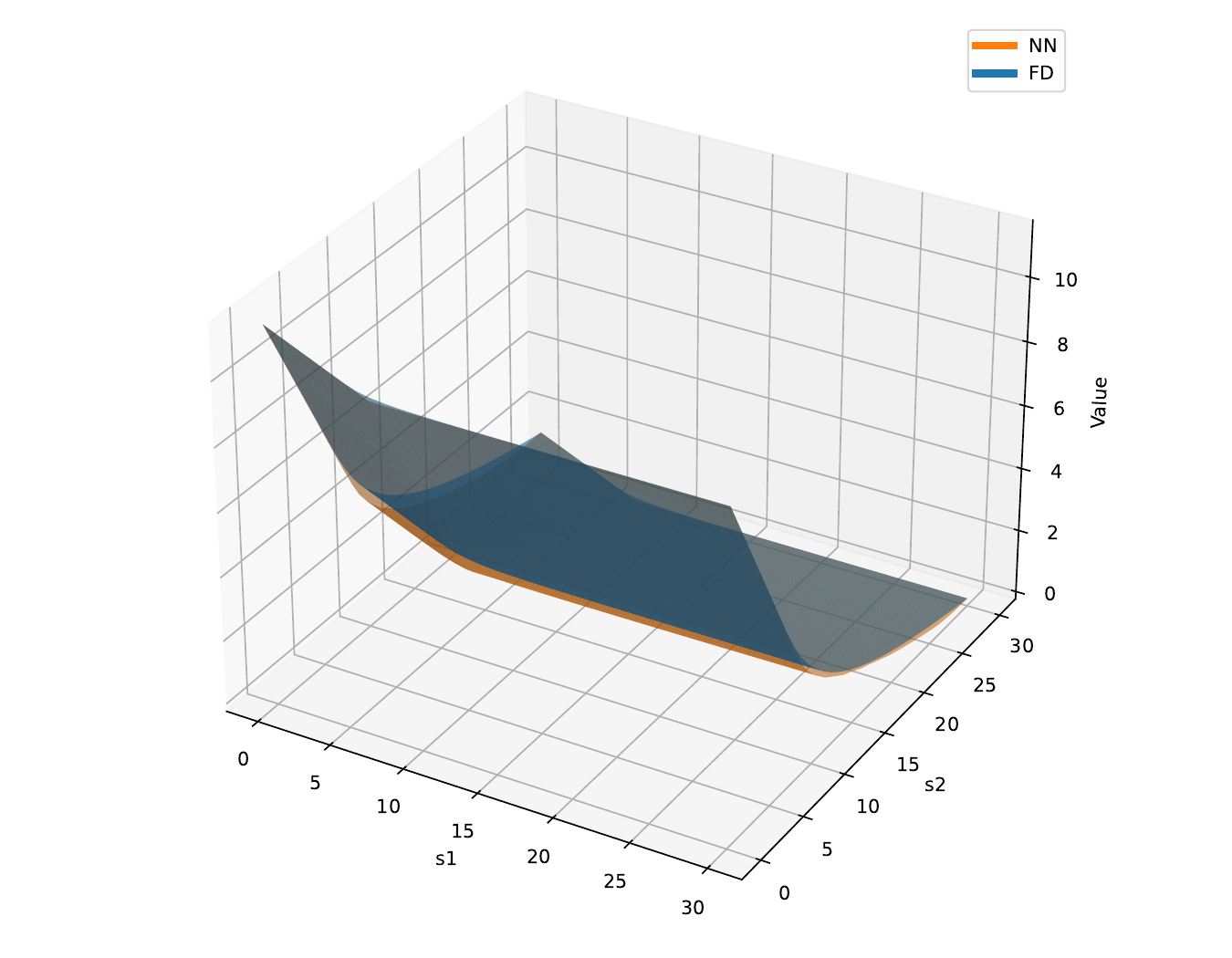}
\end{subfigure}
\caption{Multi American Put with $N=5$: Comparison of the neural network value and the finite-difference value along the diagonal $x_1=x_2=x_3=x_4=x_5$ and on the plane $x_2=x_3=x_4=x_5$ at time $t=0$ (NN stands for the neural network, payoff for the value at terminal time, and FD for finite difference).}
\label{fig:put-value}
\end{figure} 

\begin{figure}
\centering
    \includegraphics[width=.5\linewidth]{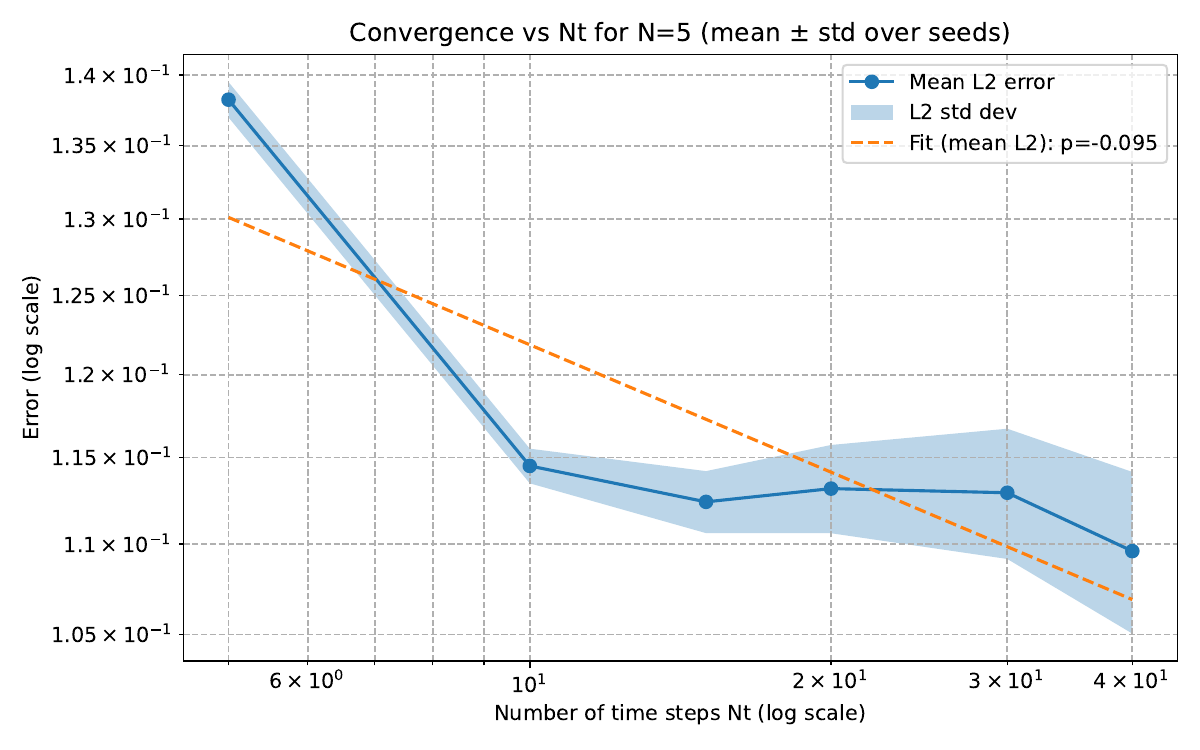}
\caption{Multi American Put with $N=5$: error as a function of the time discretization (mean $\pm$ standard deviation).}
\label{fig:put-time}
\end{figure} 

In Figure~\ref{fig:put-value}, we compare the prices of the multiple American Put ($N=5$) respectively given by our neural network and by finite differences (since the puts are independent, we compute the value of each option separately before summing them). More precisely, on the left, we plot the ``diagonal function'', that is, $x \mapsto V_0(x, \dots, x)$ where $x$ is the price of an individual stock; on the right, we plot the value on the plane $x_2 = \dots = x_N$, i.e., $(x,y) \mapsto V_0(x, y, \dots, y)$. We observe that the neural network approximates the reference solution very well by finite difference. 
In Figure~\ref{fig:put-time}, we report the $L^2$ error along the diagonal $x_1=\dots=x_N$ between the neural network prediction and the reference value, when the neural network is trained with increasing number of time steps. We observe that the error decreases. We plot the mean (full line) plus/minus standard deviation (shaded area) over 5 runs. The dashed orange line is obtained by fitting a linear model on the data points in a log-log scale. The error decreases as a function of the number of time steps, and it does so at a slower rate than the one predicted by our result (Proposition~\ref{prop:discrete-continuous2}). This is probably due to numerical errors in the optimization procedure that were not accounted for in the theoretical bound.

\paragraph*{A nonlinear example. } The following example is a sanity check to verify that our algorithm also works for ``non-separable'' utilities and to visualize the impact of the extra error added by Algorithm~\ref{alg:DBMSP}. We define:
$$ V_0(\bx) := \sup_{\boldsymbol{\t} \in \cT_{[0,T]}^N} \dbE\Big[\log\Big( 1 + \frac{1}{N}\sum_{j=1}^N S_{\t_j}^{j} \Big)\Big],$$
where the processes $S^1, \dots, S^N$ are defined as in the previous example, with the extra condition that $\m_j \le 0$ for all $j \in [N]^*$. We also introduce:
$$ V_0^h(\bx) := \sup_{\boldsymbol{\t} \in \cT_{p}^N} \dbE\Big[\log\Big( 1 + \frac{1}{N}\sum_{j=1}^N S_{\t_j}^{j,h} \Big)\Big],$$

Note that $V_0$ can be computed explicitly, as by Jensen inequality:
$$ V_0(x) \le \sup_{\boldsymbol{\t} \in \cT_{[0,T]}^N} \log\Big( \dbE\Big[1 + \frac{1}{N}\sum_{j=1}^N S_{\t_j} \Big]\Big) \le \log\Big( 1 + \frac{1}{N} \sum_{j=1}^N x_j \Big), $$
coming from the fact that the processes $\{S^j\}_{j \in [N]^*}$ are super-martingales under the condition $\m_j \le 0$. Since the above upper bound is reached for $\t_j = 0$ for all $j \in [N]^*$ (i.e., the optimal strategy is to stop immediately), we conclude that $V_0(\bx) = \log\Big( 1 + \frac{1}{N} \sum_{j=1}^N x_j \Big)$. We may proceed the same way with $V_0^h$, and we have in fact $V_0 = V_0^h$. The goal of this example is then to visualize the error due to Algorithm~\ref{alg:DBMSP}, and we therefore introduce:
$$ \tilde V_0^h(\bx) := \sup_{\boldsymbol{\t} \in \tilde \cT_{p}^N} \dbE\Big[\log\Big( 1 + \frac{1}{N}\sum_{j=1}^N S_{\t_j}^{j,h} \Big)\Big].$$

\begin{figure}
\centering
\includegraphics[scale=0.5]{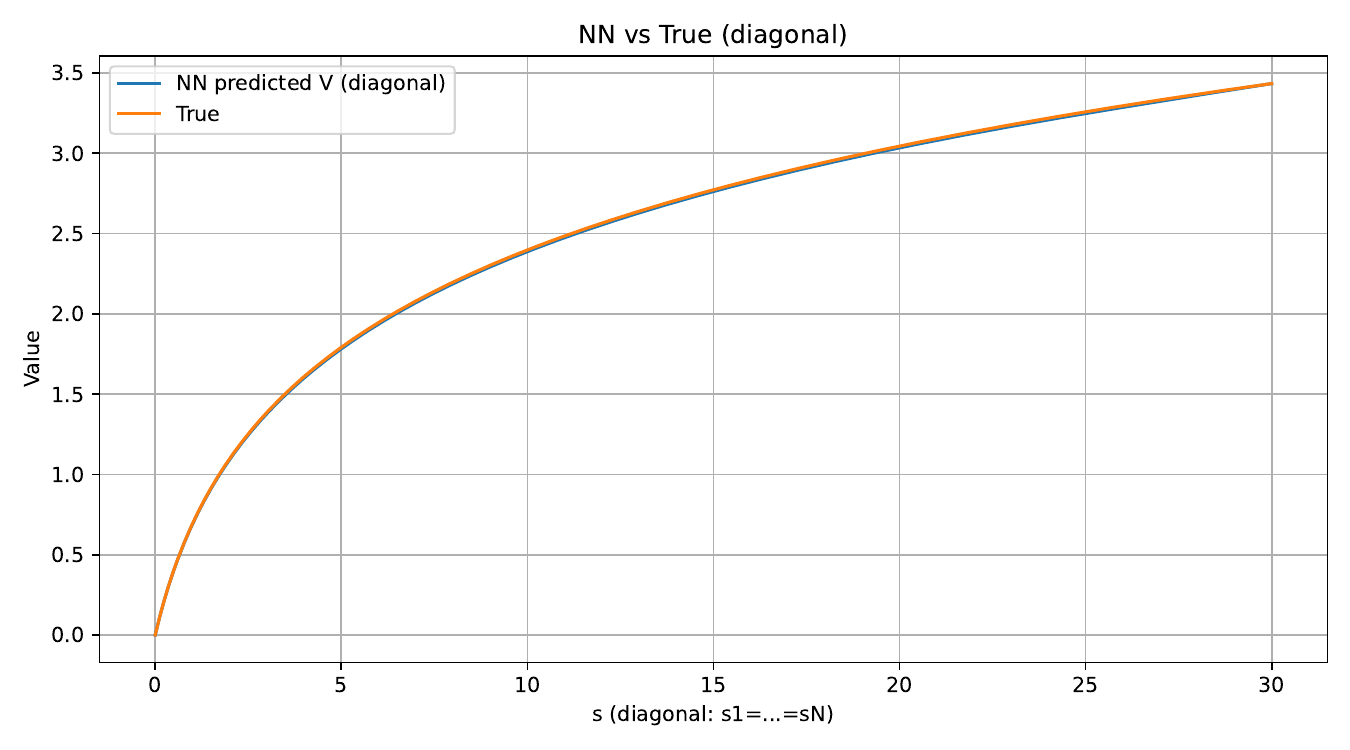}
\caption{Non-linear example with $N=5$: value functions obtained by neural network versus true value function along the diagonal $x_1 = \dots = x_N$.}
\label{fig:log-value}
\end{figure}

\begin{figure}
\centering
\includegraphics[scale=0.5]{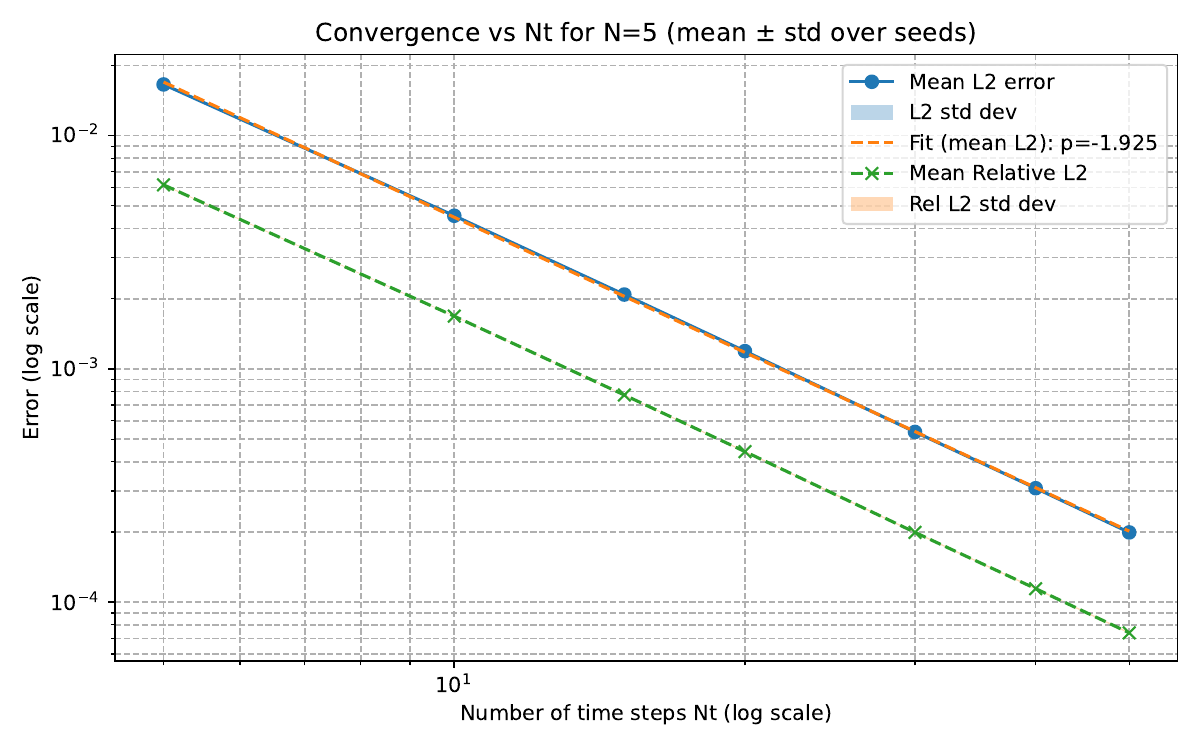}
\caption{Non-linear example with $N=5$: $L^2$ error as a function of the time discretization (mean $\pm$ standard deviation) and relative $L^2$ error.}
\label{fig:log-time}
\end{figure}

Figure~\ref{fig:log-value}, compares the actual value function and $\tilde V_0$, on the diagonal $\bx = (x, \dots, x)$ and for $N=5$.
In Figure~\ref{fig:log-time}, we report the $L^2$ error along the diagonal $x_1=\dots=x_N$ between the neural network prediction and the true value, when the neural network is trained with increasing number of time steps. We observe that the error decreases. We plot the mean (full line) plus/minus standard deviation (shaded area) over 5 runs. The orange dashed line is obtained by fitting a linear model on the data points in a log-log scale. The green dashed line shows the relative $L^2$ error computed by taking the $L^2$ norm of the difference between the values divided by the true value at each point along the diagonal (on a grid). We observe that the error decreases as a function of the number of time steps, and it does so at a faster rate than the upper bound we provided (Proposition~\ref{prop:discrete-continuous2}). This is probably due to the fact that the optimal strategy consists in stopping immediately and the value function is very smooth.

\smallskip

\noindent\textbf{Acknowledgements.} 
M.L. acknowledges that computing resources were provided by NYU Shanghai HPC.

\appendix 

\section{Technical results}\label{sec:appendix}

\begin{lemma}\label{lem:tildetau}
    Let $\bm{\t} = (\t_1, \dots, \t_N) \in \cT_p^N$, which denotes here the set of $N$-uplets of $\{t_0, \dots, t_p\}$-valued stopping times. There exists $\tilde{\bm{\t}} \in \tilde \cT_p^N$ such that $| \bm{\t} - \tilde{\bm{\t}} |_\infty \le Nh$, with $h := \min_{k \in [p-1]} |t_{k+1} - t_k |.$
\end{lemma}
\proof
Recall the following notation: $I_t^k = \1_{t \le \t_k}$ for $k \in [N]^*$ and $t \in \{ t_1, \dots, t_p \}$. For every time $t_n$, we introduce a set $\cJ_{t_n}$, which can be interpreted as the set of indices left to stop at time $t_n$. At every time, only the smallest element of the set, that is, $\min \cJ_{t_n}$, is stopped. More precisely, we define the set-valued random sequence $(\cJ_{t_n})_{n \in [p]}$ as follows: 
$$ 
    \cJ_{t_0} := \{ k \in [N]^*: I_{t_0}^k - I_{t_1}^k = 1 \},
$$
and, for $n \in [p-1]^*$,
$$ 
    \cJ_{t_ n} = \big( \cJ_{t_{n-1}} \backslash \{ \min \cJ_{t_{n-1}} \} \big) \cup \{ k \in [N]^*: I_{t_n}^k - I_{t_{n+1}}^k = 1 \}. 
$$
We also set $\cJ_p = \emptyset$. Then, for all $k \in [N]^*$, we introduce:
\bea\label{def-tildetau}
 \tilde \t_k := \min\{ t_n \ge 0: \min \cJ_{t_n} = k\} 
 \wedge t_p. 
 \eea
 We now verify that $\tilde{\bm{\t}} := \{\tilde\t_1, \dots, \tilde\t_N \}$ satisfies the desired properties. \\

First, since $I^k$ is a predictable process, by induction we see that the set-valued process  $t \mapsto \cJ_t$ is adapted to filtration $\dbF^N$. By~\eqref{def-tildetau}, it is then clear that $\tilde \t_k$ is a $\dbF^N$-stopping time. Moreover, for $l \neq k$ and $n \in [p-1]$, we have: 
\begin{align*}
    \{ \tilde \t_k = \tilde \t_l \} =& \bigcup_{n=0}^p \{ \tilde \t_k = \tilde \t_l = t_n \} 
    = \Big(\bigcup_{n=0}^{p-1} \{ k = \min \cJ_{t_n} = l \} \Big) \cup \{ \tilde \t_k = \tilde \t_l = t_p \} \\
    =& \{ \tilde \t_k = \tilde \t_l = t_p \},
\end{align*}
which implies that $\tilde{\bm{\t}} \in \tilde{\cT}_p^N$. Finally, we estimate $| \t_k - \tilde \t_k|$. Observe that:
$$ 
    \t_k = \min\{ t_n \ge t_0 : k \in \cJ_{t_n} \} \wedge t_p.
$$
Note that we also have:
$$ \tilde \t_k = \min\{ t_n \ge t_0 : k \in \cJ_{t_n} \ \mbox{and} \ k \notin \cJ_{t_{n+1}} \} \wedge t_p, $$
from which we deduce that:
\bea\label{tilde-estimate}
    | \t_k - \tilde \t_k | =  | \{ n \in [p] : k \in \cJ_{t_n} \} |h.
\eea
Now, observe that if $\cJ_{t_n}$ is nonempty, then at every time $m \ge n$, we keep removing one index from $\cJ_{t_m}$ until it is empty, although new indices may be added. However, no index $k \in [N]^*$ can return to $\cJ$ after exiting it. Therefore, since there is at most $N$ indices to pass through $\cJ$, we have, $\dbP$-almost surely:
$$ \cJ_{t_n} \neq \emptyset \Rightarrow \cJ_{(t_n + Nh) \wedge t_p } = \emptyset. $$
This implies that $| \{ n \in [p] : k \in \cJ_{t_n} \} | \le N$, and by~\eqref{tilde-estimate} we conclude the proof. \qed

\textbf{Proof of Lemma~\ref{lem:Euler-error}}
Without loss of generality, we assume that $\si_0 = 0$. Fix $\bm{\tau} \in \cT_p^N$. For all $k \in [N]^*$ and $t \in [0,T]$, we have:
\begin{align*}
    &\dbE\big[ \sup_{u \le t} | X_u^k - X_u^{k,h} |^2  \big] 
    \\
    \le& 2\dbE\Big[T \int_0^t | b_k(s, \bm{X}_s) - b_k(t_p(s), \bm{X}_{t_p(s)}^h) |^2 ds +   \int_0^t | \si_k(s, \bm{X}_s) - \si_k(t_p(s), \bm{X}_{t_p(s)}^h) |^2 ds \Big] \\
    \le& 2\dbE\Big[\int_0^t \big(T|  b_k(t_p(s), \bm{X}_{t_p(s)}) - b_k(t_p(s), \bm{X}_{t_p(s)}^h) |^2 + | \si_k(t_p(s), \bm{X}_{t_p(s)}) - \si_k(t_p(s), \bm{X}_{t_p(s)}^h) |^2\big) ds \Big] \\
    &+2\dbE\Big[\int_0^t \big(T|  b_k(s, \bm{X}_{s}) - b_k(t_p(s), \bm{X}_{t_p(s)}^h) |^2 + | \si_k(s, \bm{X}_{s}) - \si_k(t_p(s), \bm{X}_{t_p(s)}) |^2\big) ds \Big] \\
    \le& 2L(T+1)\int_0^t \dbE\big[ \sup_{u \le s} |\bX_u - \bX_u^h |^2 \big]  ds \\
    &+ 2L(T+1)\int_0^t \dbE\big[ |s - t^p(s) |^{2\beta} + |\bX_s - \bX_{t^p(s)} |^2 \big]  ds,
\end{align*}
where we successively used BDG inequality, the fact that $| I^k | \le 1$, the Lipschitz-continuity of the coefficients $b_k$ and $\si_k$ in $\bx$ and their $\beta$-Hölder-continuity in $t$. We then deduce from Gronwall's lemma that:
\begin{equation}\label{gronwall}
 \dbE\big[ \sup_{t \le T} | \bX_t - \bX_t^{h} |^2  \big] \le C_{T,N}\Big( h^{2\b} + \int_0^t \dbE\big[ |\bX_s - \bX_{t^p(s)} |^2 \big]ds  \Big). 
 \end{equation}
 Now observe that:
 \begin{align*}
\dbE\big[ |X_s^k - X_{t^p(s)}^k |^2 \big] \le& 2\dbE\Big[ \int_{t^p(s)}^s \big(T| b_k(r, \bX_r) |^2 + | \si_k(r, \bX_r) |^2 \big)dr \Big] \\
\le& h C_T \Big( 1 + \dbE\Big[ \sup_{t \le T} |\bX_t|^2 \Big] \Big) 
\le h C_{T,N}
 \end{align*}
 Using again BDG inequality and the fact that $| I^k | \le 1$, we have:
 $$ \dbE\Big[ \sup_{u \le t} |X_u^k|^2 \Big] \le C_T\int_0^t \big(1 + \dbE\big[ \sup_{u \le s} |X_u|^2 \big] \big)ds, $$
 from which we deduce by Gronwall's lemma again that the second order moment of $\bX$ are bounded independently from the stopping policy $\bm{\tau}$. Therefore:
\begin{align*}
\dbE\big[ |X_s^k - X_{t^p(s)}^k |^2 \big] 
\le h C_{T,N}
 \end{align*}  
  We finally obtain the desired result by plugging the above estimate into~\eqref{gronwall} and by observing that none of the constants involves $\bm{\tau}$.
\qed

\begin{lemma}\label{lem:continuity-tau}
Let $\bm{\tau}, \tilde{\bm{\tau}} \in \cT_{[0,T]}^N$. We denote by $\bX := (X^1, \dots, X^N)$ (resp. $\tilde \bX := (\tilde X^1, \dots, \tilde X^N)$ the dynamics~\eqref{dynamics} controlled by $\boldsymbol{\t}$ (resp. $\tilde{\boldsymbol{\t}}$). Then there exists $C_{T,N} \ge 0$ (depending on $N$ and $T$) such that:
$$ \dbE\Big[ \sup_{t \in [0,T]} | \bX_t - \tilde \bX_t |^2 \Big] \le C_{T,N} \sum_{k=1}^N \sqrt{\dbE\big[| \t^k - \tilde \t^k |\big]} $$
\end{lemma}
\proof
Without loss of generality, we write the proof for $\si_0 = 0$. Fix $t \in [0,T]$ and $k \in [N]^*$. Using convexity and Burkholder-Davis-Gundy inequality, we have:
\begin{align*}
  \dbE\Big[ \sup_{u \le t} | X_u^k - \tilde X_u^k |^2 \Big] \le C\Big( \int_0^t \dbE\Big[ | b_s^k I_s^k  - \tilde b_s^k \tilde I_s^k |^2 \Big] ds + \int_0^t \dbE\Big[ | \si_s^k I_s^k  - \tilde \si_s^k \tilde I_s^k |^2 \Big] ds\Big),
\end{align*}
where we denote $\f_s^k := \f_k(s, \bX_s)$, $\tilde \f_s^k := \f_k(s, \tilde \bX_s)$ for $\f \in \{b, \si\}$, $I_s^k := \1_{s < \t_k}$ and $\tilde I_s^k := \1_{s < \tilde \t_k}$. Now observe that:
\begin{align*}
    \dbE\Big[ | b_s^k I_s^k  - \tilde b_s^k \tilde I_s^k |^2 \Big] &\le \dbE\Big[ | b_s^k |^2 | I_s^k  - \tilde I_s^k | \Big] + \dbE\Big[ | b_s^k - \tilde b_s^k |^2 \Big] \\
    &\le C\dbE\Big[\big( 1 + | X_s^k |^2\big) | I_s^k  - \tilde I_s^k | \Big] + \dbE\Big[ | X_s^k - \tilde X_s^k |^2 \Big] \\
    &\le C\sqrt{\dbE\Big[\big( 1 + | X_s^k |^4\big)\Big]}\sqrt{\dbE\Big[ | I_s^k  - \tilde I_s^k | \Big]} + \dbE\Big[ | X_s^k - \tilde X_s^k |^2 \Big].
\end{align*}
Using estimates similar to the proof of Lemma~\ref{lem:Euler-error}, we can show that $\dbE\Big[| X_s^k |^4\Big]$ is bounded by a constant depending on $N$ and $T$ only. Then, writing the same inequalities for the term in $\si$, observing that
$$ \int_0^t \sqrt{\dbE\Big[ | I_s^k  - \tilde I_s^k | \Big]} \le \sqrt{T  \dbE\Big[ \int_0^T | I_s^k  - \tilde I_s^k |ds \Big]} =  \sqrt{T  \dbE\Big[ | \tau_k - \tilde \tau_k | \Big]}, $$
and using Gronwall's Lemma, we derive the desired result.
\qed 

\begin{lemma}\label{lem:error}
Let Assumption~\ref{assum:coef-diffusion} hold. Let $\bm{\t}, \tilde{\bm{\t}} \in \cT_p^N$ and denote by $(\bX^h, \bI)$ and $(\tilde{\bX}^h, \tilde{\bI})$ the corresponding processes defined by~\eqref{eq:euler-scheme}, with $\bX_0 = \tilde{\bX}_0$. Then we have: 
$$ 
    \dbE\Big[ \max_{n \in [p]} | \bX_n^h - \tilde{\bX}_n^h |^2\Big] \le C_{T,N} \sqrt{ \dbE\big[ | \boldsymbol{\tau} - \tilde{\bm{\tau}}| \big] } \q \mbox{for all} \ n \in [p], 
$$
where the constant $C_{T,N}$ depends on $T$ and $N$ only.
\end{lemma}
\proof
The result can be proven by adapting the same estimates as in the proof of Lemma~\ref{lem:continuity-tau} to the dynamics~\eqref{eq:euler-scheme}.
\qed

\bibliography{Bibliography}

\end{document}